\documentclass[12pt,reqno]{amsart}

\usepackage{graphicx}
\graphicspath{ {./images/} }
\usepackage{bm}

\usepackage{color}

\usepackage{amssymb}

\usepackage{amsfonts}

\usepackage{amsmath}

\usepackage{amsthm}

\usepackage[all]{xy}

\usepackage[mathlines]{lineno}

\usepackage{hyperref}

\newtheorem{theorem}{Theorem}[section]

\newtheorem{corollary}[theorem]{Corollary}

\newtheorem{lemma}[theorem]{Lemma}

\newtheorem{proposition}[theorem]{Proposition}

\newtheorem{Definition}[theorem]{Definition}

\newtheorem{Example}[theorem]{Example}

\newtheorem{Remark}[theorem]{Remark}

\newenvironment{remark}{\begin{Remark}\begin{em}}{\end{em}\end{Remark}}

\def\diag{{\mbox{diag}\,}}

\DeclareMathOperator{\tr}{tr}

\begin{document}

\title[Revisit on spectral geometric mean]{Revisit on spectral geometric mean}

\author{Luyining Gan}
\address{Department of Mathematics and Statistics, University of Nevada Reno, Reno, NV 89557, USA}
\email{lgan@unr.edu}

\author{Sejong Kim}
\address{Department of Mathematics, Chungbuk National University, Cheongju 28644, Korea}
\email{skim@chungbuk.ac.kr}

\maketitle

\begin{abstract}
In this paper we introduce the limit, unique solution of the nonlinear equations, geodesic property, tolerance relations and pinch on the spectral geometric mean for two positive definite operators. We show that the spectral geometric mean is a geodesic with respect to some semi-metric. We also prove that the tolerance relation on determinant one matrices can be characterized by the spectral geometric mean. Moreover, two positive tuples can be pinched by the spectral geometric mean.
\end{abstract}

\medskip
\noindent \textit{2020 Mathematics Subject Classification}
15A42,  
15B48,  
47A64,  
53C20.  	

\noindent \textit{Key words and phrases.} Positive definite matrix, metric geometric mean, spectral geometric mean, invariant tolerance relation.

\section{Introduction}

Averaging operations has been studied a lot in matrix theory ad operators. Many notions of means of positive definite matrices have been posed and studied extensively.
The metric geometric mean of two positive definite matrices $A$ and $B$ is defined as
\begin{equation*}
A \# B := A^{1/2}(A^{-1/2}BA^{-1/2})^{1/2}A^{1/2},
\end{equation*}
which was first introduced by Pusz and Woronowicz \cite {PW75} in 1975 and further studied by Kubo and Ando~\cite{KA79} in the 1970s. Since then, it has attracted a great deal of attention in the past several decades. Though the definition looks awkward, it is indeed a natural generalization of the classical geometric mean $\sqrt {ab}$ of two positive numbers $a, b$ \cite{Bh}. Later Kubo and Ando~\cite{KA79} established the metric geometric mean of positive definite operators, then many properties and applications including the operator inequality have been developed.

Based on the metric geometric mean, the spectral geometric mean of two positive definite $A, B$ was introduced by Fiedler and Pt\'ak \cite {FP97}
\begin{equation*}\label{eqn:sm}
A \natural B := (A^{-1}\# B)^{1/2}A(A^{-1}\# B)^{1/2}.
\end{equation*}
They named it as spectral geometric mean because the square of $A \natural B$ is similar to $AB$, which means that the eigenvalues of their spectral mean are the positive square roots of the corresponding eigenvalues of $AB$ \cite [Theorem 3.2 and Remark 3.4] {FP97}.
As the spectral geometric mean is based on the metric geometric mean form the view of the formulations, it possesses some important properties of the metric geometric mean and has been studied ~\cite{KL15, LL04, LMW21, Li12}. However, unlike the metric geometric mean, not many results have been obtained on the spectral geometric mean. Thus, this paper aims to obtain new results on the spectral geometric mean and its extension, namely, the weighted spectral geometric mean.

The structure of this article is organized as follows. We recall the fundamental properties of the weighted spectral geometric mean, and also provide a new proof of Lie-Trotter-Kato formula for the weighted spectral geometric mean in Section~\ref{Sec:Sp-Geo}. In Section~\ref{Sec:Solution} we study at the weighted spectral geometric mean from the point of view as the unique solution of the nonlinear equations. In Section~\ref{Sec:Geodesic} we prove that the weighted spectral geometric mean is a geodesic with respect to the semi-metric $d(A, B) = 2 \Vert \log (A^{-1} \# B) \Vert$, and give the counterexample of the convexity. In Section~\ref{Sec:Tolerance} we construct the equivalent conditions that $A^{-1}$ and $B$ have two kinds of tolerance relations. We also prove in Section~\ref{Sec:Pinch} that one positive $n$-tuple can be obtained by at most $n$ spectral pinches consecutively from another positive $n$-tuple if the log-majorization relation exists.

\section{Spectral geometric mean}~\label{Sec:Sp-Geo}

Let $B(\mathcal{H})$ be the Banach space of all bounded linear operators on a Hilbert space $\mathcal{H}$ with inner product $\langle \cdot, \cdot \rangle$, and let $S(\mathcal{H}) \subset B(\mathcal{H})$ be the closed subspace of all self-adjoint linear operators. We say that $A \in S(\mathcal{H})$ is positive semi-definite (positive definite) if $\langle x, Ax \rangle \geq (>) 0$ for all (nonzero, respectively) vector $x \in \mathcal{H}$. We denote as $\mathbb{P} \subset S(\mathcal{H})$ the open convex cone of all positive definite operators. The group $\textrm{GL}$ of all invertible operators transitively acts on $\mathbb{P}$ via congruence transformation. On the finite-dimensional setting $\mathcal{H} = \mathbb{C}^{m}$ we write as $\mathbb{P}_{m}$ the open convex cone of all $m \times m$ positive definite matrices.

The weighted spectral geometric mean of $A, B \in \mathbb{P}$ is a curve defined by
\begin{displaymath}
t \in \mathbb{R} \mapsto A \natural_{t} B := (A^{-1} \# B)^{t} A (A^{-1} \# B)^{t}.
\end{displaymath}
Note that $A \natural_{0} B = A$ and $A \natural_{1} B = B$. We simply write as $A \natural B = A \natural_{1/2} B$.

We list fundamental properties of the weighted spectral geometric mean.
\begin{lemma} \label{L:SG} \cite{HL07}
Let $A, B\in \mathbb{P}$ and $s, t, u \in \mathbb{R}$.
\begin{itemize}
\item[(1)] $A \natural_{t} B = A^{1-t} B^{t}$ if $A$ and $B$ commute.
\item[(2)] $(a A) \natural_{t} (b B) = a^{1-t} b^{t} (A \natural_{t} B)$ for any $a, b > 0$.
\item[(3)] $U^{*} (A \natural_{t} B) U = (U^{*} A U) \natural_{t} (U^{*} B U)$ for any unitary operator $U$.
\item[(4)] $A \natural_{t} B = B \natural_{1-t} A$.
\item[(5)] $(A \natural_{s} B) \natural_{t} (A \natural_{u} B) = A \natural_{(1-t)s + t u} B$.
\item[(6)] $(A \natural_{t} B)^{-1} = A^{-1} \natural_{t} B^{-1}$.
\end{itemize}
\end{lemma}

\begin{remark}
By \cite[Proposition 2.3]{HK} the boundedness of weighted spectral geometric mean for Loewner order has been shown:
\begin{equation} \label{E:bound}
2^{1 + t} (A + B^{-1})^{-t} - A^{-1} \leq A \natural_{t} B \leq [ 2^{1 + t} (A^{-1} + B)^{-t} - A ]^{-1}
\end{equation}
for any $A, B \in \mathbb{P}_{m}$ and $t \in (0,1)$. Using \eqref{E:bound} we provide another proof of Lie-Trotter-Kato formula for the weighted spectral geometric mean:
\begin{equation} \label{E:Lie-Trotter-Kato}
\lim_{s \to 0} \, \left( A^{s} \natural_{t} B^{s} \right)^{1/s} = \exp ((1-t) \log A + t \log B).
\end{equation}

Without loss of generality, assume that $A, B \geq \alpha I$ for some $0 < \alpha < 1$. For any $s > 0$ by \eqref{E:bound}
\begin{displaymath}
2^{1 + t} (A^{s} + B^{-s})^{-t} - A^{-s} \leq A^{s} \natural_{t} B^{s} \leq [ 2^{1 + t} (A^{-s} + B^{s})^{-t} - A^{s} ]^{-1}.
\end{displaymath}
Set $f(s) := 2^{1 + t} (A^{s} + B^{-s})^{-t} - A^{-s}$ on the interval
\begin{displaymath}
I: 0 < s < \min \left\{ t, \frac{1}{-1 + \log_{2} (\alpha^{t-1} + \alpha^{-t-1})} \right\}.
\end{displaymath}
Since $\alpha^{t-1} + \alpha^{-t-1} \geq 2 \alpha^{-1} > 2$ by the arithmetic-geometric mean inequality, we have $\log_{2} (\alpha^{t-1} + \alpha^{-t-1}) > 1$. Since $\displaystyle s < \frac{1}{-1 + \log_{2} (\alpha^{t-1} + \alpha^{-t-1})}$ and $B^{-t} \leq \alpha^{-t} I$, we have $2^{1 + \frac{1}{s}} > \alpha^{t-1} + \alpha^{-t-1}$, and
\begin{displaymath}
A^{t-1} + A^{-1/2} B^{-t} A^{-1/2} \leq A^{t-1} + \alpha^{-t} A^{-1} \leq \left( \alpha^{t-1} + \alpha^{-t-1} \right) I < 2^{1 + \frac{1}{s}} I.
\end{displaymath}
Taking congruence transformation by $A^{1/2}$ yields $A^{t} + B^{-t} < 2^{1 + \frac{1}{s}} A$.

It has been shown in \cite{BS} that the map $\displaystyle F(p) = \left( \frac{A^{p} + B^{p}}{2} \right)^{1/p}$ is operator monotone for $p \in (-\infty, -1] \cup [1, \infty)$. So for $\displaystyle \frac{t}{s} > 1$
\begin{displaymath}
\frac{A + B^{-1}}{2} \leq \left( \frac{A^{\frac{t}{s}} + B^{- \frac{t}{s}}}{2} \right)^{\frac{s}{t}}.
\end{displaymath}
Substituting $A, B$ into $A^{s}, B^{s}$, applying the monotonicity of the $t$-power map, and using the preceding argument we obtain
\begin{displaymath}
\left( \frac{A^{s} + B^{-s}}{2} \right)^{t} \leq \left( \frac{A^{t} + B^{-t}}{2} \right)^{s} < 2 A^{s}.
\end{displaymath}
Thus,
\begin{displaymath}
2^{1 + t} (A^{s} + B^{-s})^{-t} - A^{-s} = 2 \left[ \left( \frac{A^{s} + B^{-s}}{2} \right)^{-t} - \frac{1}{2} A^{-s} \right] > 0,
\end{displaymath}
that is, $f(s)$ is positive definite for any $s \in I$. Since
\begin{displaymath}
f'(s) = 2^{1+t} (-t) (A^{s} + B^{-s})^{-t-1} (A^{s} \log A - B^{-s} \log B) + A^{-s} \log A
\end{displaymath}
we can see by matrix calculus that
\begin{displaymath}
\lim_{s \to 0} \ln f(s)^{1/s} = \lim_{s \to 0} \frac{\ln f(s)}{s} = f'(0) = -t (\log A - \log B) + \log A = (1-t) \log A + t \log B.
\end{displaymath}
Hence, $\displaystyle \lim_{s \to 0} f(s)^{1/s} = \exp ((1-t) \log A + t \log B)$. Similarly, we can obtain
\begin{displaymath}
\lim_{s \to 0} [ 2^{1 + t} (A^{-s} + B^{s})^{-t} - A^{s} ]^{-1/s} = \exp((1-t) \log A + t \log B).
\end{displaymath}
Thus, we conclude \eqref{E:Lie-Trotter-Kato}.
\end{remark}

\section{Unique solution of nonlinear equations} \label{Sec:Solution}

From the Riccati Lemma: the metric geometric mean $A \# B$ is the unique positive definite solution $X \in \mathbb{P}$ of the Riccati equation $X A^{-1} X = B$, we obtain the following proposition.
\begin{proposition} \label{P:sp-mean}
Let $A, B \in \mathbb{P}$ and $t \in \mathbb{R}$. Then $A \natural_{t} B$ is a unique positive definite solution $X \in \mathbb{P}$ of the equation
\begin{displaymath}
A^{-1} \# X = (A^{-1} \# B)^{t}.
\end{displaymath}
\end{proposition}

\begin{theorem}
For non-zero $t \in \mathbb{R}$ let $G_{t}: \mathbb{P} \times \mathbb{P} \to \mathbb{P}$ be a map satisfying $G_{t}(A, A) = A$ for all $A \in \mathbb{P}$ and
\begin{equation} \label{E:condition}
G_{t}(A, B) = I \ \Longrightarrow \ B = A^{1 - \frac{1}{t}} \ \textrm{for any} \ A, B \in \mathbb{P}.
\end{equation}
Then $A \natural_{t} B$ for $A, B \in \mathbb{P}$ is a unique solution $X \in \mathbb{P}$ of the equation
\begin{displaymath}
G_{t}(A \# X^{-1}, B \# X^{-1}) = I.
\end{displaymath}
\end{theorem}

\begin{proof}
Set $U = A \# X^{-1}$ and $V = B \# X^{-1}$ for $X \in \mathbb{P}$. Then $U, V \in \mathbb{P}$ and $G_{t}(U, V) = I$, so by \eqref{E:condition}
\begin{displaymath}
V = U^{1 - \frac{1}{t}}.
\end{displaymath}
Since $X = U^{-1} A U^{-1} = V^{-1} B V^{-1}$ by the Riccati Lemma, we have
\begin{displaymath}
A = U V^{-1} B V^{-1} U = U^{\frac{1}{t}} B U^{\frac{1}{t}}.
\end{displaymath}
Again by the Riccati Lemma $U^{\frac{1}{t}} = A \# B^{-1}$, and hence, by self-duality of the metric geometric mean: $(A \# B)^{-1} = A^{-1} \# B^{-1}$
\begin{displaymath}
X = U^{-1} A U^{-1} = (A \# B^{-1})^{-t} A (A \# B^{-1})^{-t} = (A^{-1} \# B)^{t} A (A^{-1} \# B)^{t} = A \natural_{t} B.
\end{displaymath}
This completes the proof.
\end{proof}

\begin{remark}
A two-variable mean on a set $X$ is a binary operation $G: X \times X \to X$ satisfying the idempotency: $G(x, x) = x$ for all $x \in X$.  There are numerous examples of two-variable mean on the open cone $\mathbb{P}$ satisfying the idempotency and \eqref{E:condition}, such as the metric geometric mean and spectral geometric mean.
\end{remark}

The following show that a unique solution for the system of equations is given by a pair of metric geometric mean and spectral geometric mean.
\begin{theorem}
For any $A, B \in \mathbb{P}$ and $t \in \mathbb{R}$
\begin{displaymath}
\left\{
  \begin{array}{ll}
    A = X^{-t} Y X^{-t} & \hbox{$(1)$} \\
    B = X^{1-t} Y X^{1-t} & \hbox{$(2)$}
  \end{array}
\right.
\end{displaymath}
has a unique solution $(X, Y) = (A^{-1} \# B, A \natural_{t} B)$.
\end{theorem}

\begin{proof}
By (1) $Y = X^{t} A X^{t}$, and then (2) reduces to $B = X A X$. By the Riccati Lemma $X = A^{-1} \# B$, so
\begin{displaymath}
Y = (A^{-1} \# B)^{t} A (A^{-1} \# B)^{t} = A \natural_{t} B.
\end{displaymath}
\end{proof}

\section{Geodesic property}~\label{Sec:Geodesic}

A curve $\gamma: I \to X$ from an interval $I$ to the metric space $(X, d)$ is called a \emph{geodesic} if there exists a constant $c \geq 0$ such that for any $t \in I$ there is a neighborhood $J$ containing $t$ such that
\begin{displaymath}
d(\gamma(t_{1}), \gamma(t_{2})) = c |t_{1} - t_{2}|
\end{displaymath}
for any $t_{1}, t_{2} \in J$. This is a generalized notion of the geodesic in Riemannian manifold. If the above equality with $c = 1$ holds for any $t_{1}, t_{2} \in I$ the geodesic is called a minimizing geodesic or shortest path. The arc length of a minimizing geodesic between two points defines a distance function for any Riemannian manifold, which makes it into a metric space.

We consider the map $d: \mathbb{P} \times \mathbb{P} \to [0, \infty)$ given by
\begin{displaymath}
d(A, B) = 2 \Vert \log (A^{-1} \# B) \Vert,
\end{displaymath}
where $\Vert \cdot \Vert$ denotes the operator norm. From \cite[Lemma 2]{Kim21} it is a semi-metric: it satisfies the axioms of metric except the triangle inequality. Furthermore, it is invariant under the homogeneity, inversion, and unitary congruence transformation.

\begin{theorem}
The weighted geometric mean is a geodesic for the semi-metric $d$, that is,
for any $s, t \in \mathbb{R}$
\begin{displaymath}
d(A \natural_{s} B, A \natural_{t} B) = |s-t| d(A, B).
\end{displaymath}
\end{theorem}

\begin{proof}
When $s = 1$
\begin{displaymath}
\begin{split}
d(A \natural_{1} B, A \natural_{t} B) & = d(B, B \natural_{1-t} A) = 2 \Vert \log (B^{-1} \# (B \natural_{1-t} A)) \Vert \\
& = 2 \Vert \log (B^{-1} \# A)^{1-t} \Vert = |1-t| d(A, B).
\end{split}
\end{displaymath}
The first equality follows from Lemma \ref{L:SG} (4), and the third equality follows from Proposition \ref{P:sp-mean}.

Let $s \neq 1$. Then one can write $\displaystyle A \natural_{t} B = (A \natural_{s} B) \natural_{\frac{t-s}{1-s}} B$ by Lemma \ref{L:SG} (5). Then
\begin{displaymath}
\begin{split}
d(A \natural_{s} B, A \natural_{t} B) & = d \left( A \natural_{s} B, (A \natural_{s} B) \natural_{\frac{t-s}{1-s}} B \right) \\
& = 2 \left\Vert \log \left[ (A \natural_{s} B)^{-1} \# B \right]^{\frac{t-s}{1-s}} \right\Vert \\
& = 2 \left| \frac{t-s}{1-s} \right| \cdot \left\Vert \log \left[ (A \natural_{s} B)^{-1} \# B \right] \right\Vert \\
& = \left| \frac{t-s}{1-s} \right| d(A \natural_{s} B, B) \\
& = \left| \frac{t-s}{1-s} \right| \cdot |1-s| d(A, B) = |t-s| d(A, B).
\end{split}
\end{displaymath}
The second equality follows from Proposition \ref{P:sp-mean}, and the fifth equality follows from the preceding argument.
\end{proof}

\begin{remark}
The weighted metric geometric mean satisfies the convexity for the Thompson metric:
\begin{displaymath}
d_{T}(A \#_{t} B, A \#_{t} C) \leq t d_{T}(B, C)
\end{displaymath}
for $A, B, C \in \mathbb{P}$ and $t \in [0,1]$, where $d_{T}(A, B) = \Vert \log A^{-1/2} B A^{-1/2} \Vert$ denotes the Thompson metric. It gives a generalized convexity
\begin{displaymath}
d_{T}(A \#_{s} B, C \#_{t} D) \leq (1-s) d_{T}(A, C) + s d_{T}(B, D) + |s-t| d_{T}(C, D)
\end{displaymath}
for $s, t \in [0,1]$, and plays a very important role to extension of two-variable geometric mean to multi-variable geometric means as Ando-Li-Mathias mean \cite{ALM}, Bini-Meini-Poloni mean \cite{BMP}, and Cartan mean \cite{Ho, LP14}.

Unfortunately, the weighted spectral geometric mean does not satisfy the convexity for the semi-metric $d$. In other words,
\begin{displaymath}
d(A \natural_{t} B, A \natural_{t} C) \leq t d(B, C)
\end{displaymath}
does not always hold for $A, B, C \in \mathbb{P}$ and $t \in [0,1]$. For a counterexample, let
\begin{displaymath}
A =
\left(
  \begin{array}{cc}
    12.9638 & 8.0820 \\
    8.0820 & 10.9249 \\
  \end{array}
\right), \
B =
\left(
  \begin{array}{cc}
    11.3531 & 9.1847 \\
    9.1847 & 11.9930 \\
  \end{array}
\right), \
C =
\left(
  \begin{array}{cc}
    21.8929 & -10.7568 \\
    -10.7568 & 39.9958 \\
  \end{array}
\right).
\end{displaymath}
Using MATLAB we can get $d(A \natural_{1/3} B, A \natural_{1/3} C) \approx 0.9328$ and $\frac{1}{3} d(B, C) \approx 0.9266$.
\end{remark}

For the uniqueness of a minimal geodesic for the Thompson metric on $\mathbb{P}_{2}$ we have from \cite[Lemma 2.4]{CGL}
\begin{equation} \label{E:linear-sum}
A \#_{t} B = L_{1-t} (\lambda) A + L_{t} (\lambda) B
\end{equation}
for any $A, B \in \mathbb{P}_{2}$ and $t \in [0,1]$, where $\lambda$ is an eigenvalue of $A B^{-1}$ and
\begin{displaymath}
L_{t}(\lambda) =
\left\{
  \begin{array}{ll}
    \displaystyle \frac{\lambda^{t} - \lambda^{-t}}{\lambda - \lambda^{-1}}, & \hbox{$\lambda \neq 1$;} \\
    t, & \hbox{$\lambda = 1$.}
  \end{array}
\right.
\end{displaymath}
The formula \eqref{E:linear-sum} is a generalization of the consequence in \cite[Proposition 4.1.12]{Bh}
\begin{displaymath}
A \# B = \frac{A + B}{\sqrt{\det (A + B)}}
\end{displaymath}
for any $A, B \in \mathbb{P}_{2}$ with determinant $1$. It is a natural question whether or not the spectral geometric mean on $\mathbb{P}$ can be written as a linear sum. It is obvious for commuting $A, B \in \mathbb{P}_{2}$ that the spectral geometric mean has a linear form:
\begin{displaymath}
A \natural_{t} B = A^{1-t} B^{t} = A \#_{t} B = L_{1-t} (\lambda) A + L_{t} (\lambda) B.
\end{displaymath}
On the other hand, it does not hold in general for non-commuting positive definite operators $A$ and $B$.

\begin{theorem}
The linear sum for two-variable spectral geometric mean $A \natural_{t} B$ of non-commuting $A, B \in \mathbb{P}$ for $t \in (0,1)$ is not available.
\end{theorem}

\begin{proof}
For $A, B \in \mathbb{P}$ which $A B \neq B A$, assume that there exist $x, y > 0$ such that $A \natural_{t} B = x A + y B$ for $t \in (0,1)$. Then by Proposition \ref{P:sp-mean}
\begin{displaymath}
\begin{split}
(A^{-1} \# B)^{t} & = A^{-1} \# (x A + y B) \\
& = A^{-1/2} (x A^{2} + y A^{1/2} B A^{1/2})^{1/2} A^{-1/2}.
\end{split}
\end{displaymath}
Since the map $A \mapsto A^{1/2}$ is operator concave on $\mathbb{P}$,
\begin{displaymath}
(A^{-1} \# B)^{t} \geq A^{-1/2} \left[ x A + y (A^{1/2} B A^{1/2})^{1/2} \right] A^{-1/2} = x I + y A^{-1} \# B.
\end{displaymath}
It says that all positive eigenvalues $\lambda$ of a positive definite operator $A^{-1} \# B$ must satisfy
\begin{displaymath}
\lambda^{t} \geq x + y \lambda.
\end{displaymath}
Set $f(\lambda) := \lambda^{t} - x - y \lambda$. Then $f''(\lambda) = t(t-1) \lambda^{t-2} < 0$, and $f(\lambda) \to -x$ as $\lambda \to 0^{+}$. Since $\lambda$ is arbitrarily positive, it does not hold that $f(\lambda) \geq 0$.
\end{proof}

\section{Tolerance relation}~\label{Sec:Tolerance}

A tolerance relation on a set, first recognized by Poincar\'{e} \cite{Poin}, is a reflexive and symmetric relation. One can see that it is like a congruence, except the assumption of transitivity. Tolerance relations provide a convenient tool for studying indiscernibility and indistinguishability phenomena. We study in this section the effect of tolerance relation on spectral geometric means.

For $A, B \in \mathbb{P}_{m}$ we define $A \sim B$ if there exist $a, b > 0$ such that
\begin{center}
$\sqrt{a b} = \det(A^{-1} B)^{\frac{1}{m}}$ \ and \ $\sigma(A^{-1} B) = \{ a, b \}$.
\end{center}
Note that $\sim$ is a tolerance relation on $\mathbb{P}_{m}$. Furthermore, it is invariant under Riemannian isometries such as the inversion and congruence transformation.

\begin{remark} \label{R:tolerance}
If $A$ and $B$ are linearly independent and $A \sim B$, then $m$ is even, and $A \sim B$ if and only if $A^{-1/2} B A^{-1/2}$ is similar to
\begin{displaymath}
D_{m} (a, b) := \textrm{diag}(\underset{m/2}{\underbrace{a, \dots, a}}, \underset{m/2}{\underbrace{b, \dots, b}})
\end{displaymath}
for some $a, b > 0$.
\end{remark}

From the known result, we obtain
\begin{lemma} \cite{Lim21} \label{L:supplement}
For $A, B \in \mathbb{P}_{m}$ with $\det A = \det B = 1$, $A \sim B$ if and only if
\begin{displaymath}
A \# B = \frac{1}{\sqrt[m]{\det(A + B)}} (A + B).
\end{displaymath}
\end{lemma}
We provide several equivalent conditions that $A^{-1}$ and $B$ have the tolerance relation.

\begin{theorem} \label{T:tolerance}
Let $A, B \in \mathbb{P}_{m}$ with $\det A = \det B = 1$, and $t \in \mathbb{R}$. Then the following are equivalent.
\begin{itemize}
\item[(1)] $A^{-1} \sim B$;
\item[(2)] $\displaystyle A^{-1} \# B = \frac{1}{\sqrt[m]{\det(I + A B)}} (A^{-1} + B)$;
\item[(3)] $\displaystyle A \natural_{t} B = \frac{1}{\sqrt[m]{\det(I + A B)^{2t}}} (A^{-1} + B)^{t} A (A^{-1} + B)^{t}$.
\end{itemize}
\end{theorem}

\begin{proof}
From Lemma \ref{L:supplement}, it is obvious that (1) and (2) are equivalent.
\begin{itemize}
\item[(2)] $\Rightarrow (3)$: Applying (2) to definition of the weighted spectral geometric mean, we can easily obtain the formula (3).

\item[(3)] $\Rightarrow (2)$: If $t = 0$ it is trivial. Assume that (3) holds for any non-zero real number $t$. Set $X := (A^{-1} \# B)^{t}$. Then by definition of the weighted spectral geometric mean
\begin{displaymath}
X A X = \frac{1}{\det(I + A B)^{2t/m}} (A^{-1} + B)^{t} A (A^{-1} + B)^{t}.
\end{displaymath}
By the Riccati Lemma and using the definition of metric geometric mean,
\begin{displaymath}
\begin{split}
X & = A^{-1} \# \left[ \frac{1}{\det(I + A B)^{2t/m}} (A^{-1} + B)^{t} A (A^{-1} + B)^{t} \right] \\
& = \frac{1}{\det(I + A B)^{t/m}} (A^{-1} + B)^{t}.
\end{split}
\end{displaymath}
Taking the $1/t$ power on both sides, we obtain (2).
\end{itemize}
\end{proof}

\begin{corollary} \label{C:sp-mean}
Let $A, B \in \mathbb{P}_{m}$ with $\det A = \alpha$ and $\det B = \beta$, and $t \in (0,1)$. Then $A^{-1} \sim B$ implies
\begin{displaymath}
A \natural_{t} B = \frac{(\alpha \beta)^{3t/m}}{\det((\alpha \beta)^{1/m} I + A B)^{2t/m}} (A^{-1} + (\alpha \beta)^{-1/m} B)^{t} A (A^{-1} + (\alpha \beta)^{-1/m} B)^{t}.
\end{displaymath}
\end{corollary}

\begin{proof}
Assume that $A^{-1} \sim B$ for $A, B \in \mathbb{P}_{m}$. That is, there exist $a, b > 0$ such that
\begin{center}
$\sqrt{a b} = \det(A B)^{\frac{1}{m}}$ \ and \ $\sigma(A B) = \{ a, b \}$.
\end{center}
Let $A_{1} = \alpha^{-1/m} A$ and $B_{1} = \beta^{-1/m} B$. Then
\begin{displaymath}
\begin{split}
\sigma(A_{1} B_{1}) & = (\alpha \beta)^{-1/m} \sigma(A B) = \left\{ (\alpha \beta)^{-1/m} a, (\alpha \beta)^{-1/m} b \right\}, \\
\det(A_{1} B_{1})^{1/m} & = (\alpha \beta)^{-1/m} \det (A B)^{1/m} = (\alpha \beta)^{-1/m} \sqrt{a b}.
\end{split}
\end{displaymath}
So $A_{1}^{-1} \sim B_{1}$. Since $A_{1}, B_{1} \in \mathbb{P}_{m}$ have determinant $1$, Theorem \ref{T:tolerance} yields
\begin{displaymath}
A_{1} \natural_{t} B_{1} = \frac{1}{\sqrt[m]{\det(I + A_{1} B_{1})^{2t}}} (A_{1}^{-1} + B_{1})^{t} A_{1} (A_{1}^{-1} + B_{1})^{t}.
\end{displaymath}
By a simple calculation together with the joint homogeneity of weighted spectral geometric mean in Lemma \ref{L:SG} (2), we can obtain the desired formula.
\end{proof}

\begin{remark}
Since any $A, B \in \mathbb{P}_{2}$ have tolerance relation $A^{-1} \sim B$,
\begin{displaymath}
A \natural_{t} B = \frac{(\alpha \beta)^{t/2}}{\det((\alpha \beta)^{1/2} I + A B)^{t}} (\sqrt{\alpha \beta} A^{-1} + B)^{t} A (\sqrt{\alpha \beta} A^{-1} + B)^{t}.
\end{displaymath}
Since $\det((\alpha \beta)^{1/2} I + A B) = \sqrt{\alpha \beta} (2 \sqrt{\alpha \beta} + \tr (A B))$, we have
\begin{displaymath}
A \natural_{t} B = \left[ \frac{1}{2 \sqrt{\alpha \beta} + \tr (A B)} \right]^{t} (\sqrt{\alpha \beta} A^{-1} + B)^{t} A (\sqrt{\alpha \beta} A^{-1} + B)^{t},
\end{displaymath}
which is already provided in \cite[Corollary 2]{Kim21}. So Corollary \ref{C:sp-mean} is a generalization of \cite[Corollary 2]{Kim21} under tolerance relation.
\end{remark}

There is another invariant tolerance relation on $\mathbb{P}_{m}$. For $A, B\in \mathbb{P}_{m}$, we define $A\sigma B$ if there exist $a, b>0$ such that $\sigma (A^{-1}B) = \{ a, b\}$. From the recent result, we know that it can be characterized by the linear sum for the two-variable weighted metric geometric mean $A\#_t B$. Note that for $t\in \mathbb{R}$ and $a, b>0$, define
\[L_{a, b}(t) = \frac{ab^t - ba^t}{a-b}.
\]

\begin{lemma}\cite{Lim21} \label{L:linear}
For $A, B \in \mathbb{P}_m$ and $0<t<1$, $A\sigma B$ if and only if $t\mapsto A\#_t B$ lies in the positive linear span of $A$ and $B$, that is,
\[A\#_t B = L_{a, b}(t) A+ L_{a^{-1}, b^{-1}}(1-t) B,\]
where $\sigma(A^{-1}B) = \{ a, b\}$.
\end{lemma}

Motivated by Lemma~\ref{L:linear}, we provide several equivalent conditions that $A^{-1}$ and $B$ have the tolerance relation.

\begin{theorem}
Let $A, B \in \mathbb{P}_{m}$ and $0<t<1$. Then the following are equivalent.
\begin{itemize}
\item[(1)] $A^{-1} \sigma B$;
\item[(2)] $\displaystyle A^{-1} \# B = \frac{\sqrt{a}-\sqrt{b}}{a-b} (\sqrt{ab}A^{-1} + B)$;
\item[(3)] $\displaystyle A \natural_{t} B = \frac{\sqrt{a}-\sqrt{b}}{a-b} (\sqrt{ab}A^{-1} + B)^{t} A (\sqrt{ab}A^{-1} + B)^{t}$.
\end{itemize}
\end{theorem}

\begin{proof}
By Lemma \ref{L:linear}, it is obvious that (1) and (2) are equivalent for the case $t=1/2$.
\begin{itemize}
\item[(2)] $\Rightarrow (3)$: Applying (2) to definition of the weighted spectral geometric mean, we can easily obtain the formula (3).

\item[(3)] $\Rightarrow (2)$: If $t = 0$ it is trivial. Assume that (3) holds for any non-zero real number $t$. Set $X := (A^{-1} \# B)^{t}$. Then by definition of the weighted spectral geometric mean
\begin{displaymath}
X A X = \left(\frac{\sqrt{a}-\sqrt{b}}{a-b}\right)^{2t} (\sqrt{ab}A^{-1} + B)^t A (\sqrt{ab}A^{-1} + B)^t.
\end{displaymath}
By the Riccati Lemma and using the definition of metric geometric mean,
\begin{displaymath}
\begin{split}
X & = A^{-1} \# \left[ \left(\frac{\sqrt{a}-\sqrt{b}}{a-b}\right)^{2t}(\sqrt{ab}A^{-1} + B)^t A (\sqrt{ab}A^{-1} + B)^t \right] \\
& = \left(\frac{\sqrt{a}-\sqrt{b}}{a-b}\right)^{t} (\sqrt{ab}A^{-1} + B)^{t}.
\end{split}
\end{displaymath}
Taking the $1/t$ power on both sides, we obtain (2).
\end{itemize}
\end{proof}

\section{Spectral pinch}~\label{Sec:Pinch}

Let $\alpha, \beta \in \mathbb R^m$. We say that $\beta$ is a \emph{pinch} of $\alpha$ (see \cite{MOA11}, page 17), if
\begin{displaymath}
\beta = t \alpha +(1-t) Q \alpha
\end{displaymath}
for some $t \in [0,1]$, where $Q$ is the $m \times m$ permutation matrix that interchanges two coordinates. This means that $\beta$ is obtained by $\alpha$ under the weighted arithmetic mean $t I + (1-t)Q$. It is well known that if $\beta$ is majorized by $\alpha$, denoted by $\beta \prec \alpha$, then $\beta$ can be obtained by applying at most $m$ pinches consecutively, starting from $\alpha$. The converse is clearly true.

Later Dinh, Ahsani and Tam~\cite{DAT16} provided the definition of geometric pinch.
Let $\alpha = (\alpha_{1}, \dots, \alpha_{m}), \beta = (\beta_{1}, \dots, \beta_{m}) \in \mathbb R_+^m$, where $\mathbb R_+^m$ denotes the set of all positive $m$-tuples.
We say that $\beta$ is a \emph {geometric pinch} of $\alpha$ if
\[
\diag (\beta_1, \dots, \beta_m) = (Q^T \diag (\alpha_1, \dots, \alpha_m) Q) \#_t \diag (\alpha_1, \dots, \alpha_m)
\]
for some $t\in [0,1]$ and some transposition matrix $Q$.
They also proved the following theorem related to the log-majorization.

\begin{theorem} \cite{DAT16} \label{T:geo-pinch}
Let $\alpha, \beta \in \mathbb R_+^m$.
If $\beta$ is log-majorized by $\alpha$, denoted by $\beta \prec_{\log} \alpha$, then $\beta$ can be obtained by applying at most $m$ geometric pinches consecutively, starting from $\alpha$.
\end{theorem}

Naturally we define and consider the spectral pinch.
We say that $\beta$ is a \emph {spectral pinch} of $\alpha$ for $\alpha, \beta \in \mathbb R_+^m$ if
\[
\diag (\beta_1, \dots, \beta_m) = (Q^T \diag (\alpha_1, \dots, \alpha_m) Q)\natural_t \diag (\alpha_1, \dots, \alpha_m)
\]
for some $t \in [0,1]$ and some transposition matrix $Q$.

\begin{theorem} \label{T:sp-pinch}
Let $\alpha, \beta \in \mathbb R_+^m$.
If $\beta \prec_{\log} \alpha$, then $\beta$ can be obtained by applying at most $m$ spectral pinches consecutively, starting from $\alpha$.
\end{theorem}

\begin{proof}
	If $\beta \prec_{\log} \alpha$, that is, $\log \beta \prec \log \alpha$, then we know $\log \beta$ can be obtained by applying at most $n$ pinches consecutively starting from $\log \alpha$.
	Let $\log \hat{\alpha}$ be a pinch of $\log \alpha$. Without loss of generality, we may assume that the pinch occurs on the first two coordinates.
	So $(\hat{\alpha_1}, \hat{\alpha_2})\prec_{\log} (\alpha_1, \alpha_2)$ and thus $\hat{\alpha_1} = \alpha_1^t\alpha_2^{1-t}$ and $\hat{\alpha_2} = \alpha_2^t\alpha_1^{1-t}$ for some $t\in [0,1]$. Let $P$ denote the matrix corresponding to the transposition switching the first two coordinates. Then by Lemma \ref{L:SG} (1)
\begin{displaymath}
\begin{split}
& (P^T \diag(\alpha_1, \alpha_2, \alpha_3, \dots, \alpha_m) P) \natural_t \diag(\alpha_1, \alpha_2, \alpha_3, \dots, \alpha_m) \\
& = \diag(\alpha_2, \alpha_1, \alpha_3, \dots, \alpha_m) \natural_t \diag(\alpha_1, \alpha_2, \alpha_3, \dots, \alpha_m) \\
& = \diag(\alpha_1^t \alpha_2^{1-t}, \alpha_2^t \alpha_1^{1-t}, \alpha_3, \dots, \alpha_m).
\end{split}
\end{displaymath}
Then repeat the process to conclude that there exist $t_1, \dots, t_k \in [0,1]$ and transposition matrices $P_1, \dots, P_k$ such that
\[
\diag \alpha^{(i+1)} := (P_i^T (\diag \alpha^{(i)}) P_i) \natural_{t_i} \diag \alpha^{(i)}, \quad i=1, \dots, k,
\]
where $\alpha^{(1)} := \alpha$ and $\alpha^{(k+1)} := \beta$.
\end{proof}

\begin{remark}
Let $A, B \in \mathbb{P}_{m}$. Note that $A \prec_{\log} B$ if $\lambda(A) \prec_{\log} \lambda(B)$, where $\lambda(A)$ denotes the $m$-tuple of positive eigenvalues of $A$. We can generalize the definition of geometric and spectral pinches such as $B$ is a geometric pinch (spectral pinch) of $A$ if
\begin{center}
$A = Q^{T} B Q \#_{t} B$ \ ($A = Q^{T} B Q \natural_{t} B$ respectively)
\end{center}
for some $t \in [0,1]$ and some transposition matrix $Q$. If $A, B$ are diagonal matrices then they are the same as the original definitions so Theorem \ref{T:geo-pinch} and Theorem \ref{T:sp-pinch} hold for $A \prec_{\log} B$. Since $B$ and $Q^{T} B Q$ do not commute in general, geometric pinch and spectral pinch are different. It would be an interesting problem to find relationships between the log-majorization and geometric (spectral) pinch of $A$ and $B$.
\end{remark}

\vspace{1cm}

\textbf{Acknowledgement}

The work of L. Gan was supported by the AMS-Simons Travel Grant 2022-2024.
The work of S. Kim was supported by the National Research Foundation of Korea grant funded by the Korea government (MSIT) (No. NRF-2022R1A2C4001306).

\bibliographystyle{abbrv}
\bibliography{refs}

\begin{thebibliography}{10}

\bibitem{ALM}
T.~Ando, C.-K. Li, and R.~Mathias.
\newblock Geometric means.
\newblock {\em Linear Algebra Appl.}, 385:305--334, 2004.

\bibitem{BS}
K.~V. Bhagwat and R.~Subramanian.
\newblock Inequalities between means of positive operators.
\newblock {\em Math. Proc. Cambridge Philos. Soc.}, 83(3):393--401, 1978.

\bibitem{Bh}
R.~Bhatia.
\newblock {\em Positive definite matrices}.
\newblock Princeton Series in Applied Mathematics. Princeton University Press,
  Princeton, NJ, 2007.
\newblock [2015] paperback edition of the 2007 original [ MR2284176].

\bibitem{BMP}
D.~A. Bini, B.~Meini, and F.~Poloni.
\newblock An effective matrix geometric mean satisfying the
  {A}ndo-{L}i-{M}athias properties.
\newblock {\em Math. Comp.}, 79(269):437--452, 2010.

\bibitem{CGL}
H.~Choi, E.~Ghiglioni, and Y.~Lim.
\newblock The {K}archer mean of three variables and quadric surfaces.
\newblock {\em J. Math. Anal. Appl.}, 490(2):124321, 21, 2020.

\bibitem{DAT16}
T.~H. Dinh, S.~Ahsani, and T.-Y. Tam.
\newblock Geometry and inequalities of geometric mean.
\newblock {\em Czechoslovak Math. J.}, 66(141)(3):777--792, 2016.

\bibitem{FP97}
M.~Fiedler and V.~Pt\'{a}k.
\newblock A new positive definite geometric mean of two positive definite
  matrices.
\newblock {\em Linear Algebra Appl.}, 251:1--20, 1997.

\bibitem{Ho}
J.~Holbrook.
\newblock No dice: a deterministic approach to the {C}artan centroid.
\newblock {\em J. Ramanujan Math. Soc.}, 27(4):509--521, 2012.

\bibitem{Kim21}
S.~Kim.
\newblock Operator inequalities and gyrolines of the weighted geometric means.
\newblock {\em Math. Inequal. Appl.}, 24(2):491--514, 2021.

\bibitem{KL15}
S.~Kim and H.~Lee.
\newblock Relative operator entropy related with the spectral geometric mean.
\newblock {\em Anal. Math. Phys.}, 5(3):233--240, 2015.

\bibitem{HK}
S.~Kim and H.~Lee.
\newblock Relative operator entropy related with the spectral geometric mean.
\newblock {\em Anal. Math. Phys.}, 5(3):233--240, 2015.

\bibitem{KA79}
F.~Kubo and T.~Ando.
\newblock Means of positive linear operators.
\newblock {\em Math. Ann.}, 246(3):205--224, 1979/80.

\bibitem{LL04}
J.~Lawson and Y.~Lim.
\newblock Means on dyadic symmetric sets and polar decompositions.
\newblock {\em Abh. Math. Sem. Univ. Hamburg}, 74:135--150, 2004.

\bibitem{HL07}
H.~Lee and Y.~Lim.
\newblock Metric and spectral geometric means on symmetric cones.
\newblock {\em Kyungpook Math. J.}, 47(1):133--150, 2007.

\bibitem{LMW21}
L.~Li, L.~Moln\'{a}r, and L.~Wang.
\newblock On preservers related to the spectral geometric mean.
\newblock {\em Linear Algebra Appl.}, 610:647--672, 2021.

\bibitem{Li12}
Y.~Lim.
\newblock Factorizations and geometric means of positive definite matrices.
\newblock {\em Linear Algebra Appl.}, 437(9):2159--2172, 2012.

\bibitem{Lim21}
Y.~Lim.
\newblock Invariant tolerance relations on positive definite matrices.
\newblock {\em Linear Algebra Appl.}, 619:1--11, 2021.

\bibitem{LP14}
Y.~Lim and M.~P\'{a}lfia.
\newblock Weighted deterministic walks for the squares mean on {H}adamard
  spaces.
\newblock {\em Bull. Lond. Math. Soc.}, 46(3):561--570, 2014.

\bibitem{MOA11}
A.~W. Marshall, I.~Olkin, and B.~C. Arnold.
\newblock {\em Inequalities: theory of majorization and its applications}.
\newblock Springer Series in Statistics. Springer, New York, second edition,
  2011.

\bibitem{Poin}
H.~Poincar\'{e}.
\newblock {\em Science and Hypothesis}.
\newblock The Walter Scott Publishing Co., LTD.
\newblock with a preface by Larmor J.

\bibitem{PW75}
W.~Pusz and S.~L. Woronowicz.
\newblock Functional calculus for sesquilinear forms and the purification map.
\newblock {\em Rep. Mathematical Phys.}, 8(2):159--170, 1975.

\end{thebibliography}

\end{document}